\documentclass[12pt]{amsart}
\newcommand{\bdes}{\begin{description}}
\newcommand{\edes}{\end{description}}
\newcommand{\beqn}{\begin{equation}}
\newcommand{\eeqn}{\end{equation}}

\newcommand{\PP }{{\mathbb P}}
\newcommand{\QQ }{{\mathbb Q}}
\newcommand{\CC }{{\mathbb C}}

\newcommand{\ZZ }{{\mathbb Z}}

\newcommand{\TT }{{\mathbb T}}
\newcommand{\hb }{{\hbar}}
\newcommand{\cT}{{\mathcal T}}
\newcommand{\co}{{\mathcal O}}

\newcommand{\pl}{{\tilde {\mathcal P}}}
\newcommand{\cf}{{\text{ft}}}

\newtheorem{proposition}{Proposition}[subsection]
\newtheorem{lemma}{Lemma}[subsection]

\newtheorem{conjecture}{Conjecture}[subsection]
\newtheorem{def/th}{Definition/Theorem}[subsection]

\begin{document}
\title{Mirror symmetry and quantum cohomology of projective bundles}
\author{Artur Elezi}
\begin{abstract}
In \cite{[E]} we conjectured a relation between the quantum
$\mathcal D$-modules of a smooth variety $X$ and the
projectivisation of a direct sum of line bundles over it. In this
paper we prove the conjecture when $X$ is a semiample complete
intersection in a toric variety. We use the conjecture to show that
the relations of the small quantum cohomology ring of $X$ that come
from differential operators lift to the projective bundle. The basic
cohomology relation of the projective bundle deforms to a relation
in the small quantum cohomology.
\end{abstract}

\maketitle

\section{Introduction}

\noindent Let $Y$ be a projective manifold. Denote by $Y_{k,\beta}$
the moduli stack of rational stable maps of class $\beta\in
H_2(Y,\ZZ)$ with $k$-markings \cite {[FP]} and $[Y_{k,\beta}]$ its
virtual fundamental class \cite{[BF]},\cite{[LT]}. Throughout this
paper we will be interested mainly in $k=1$. Recall the following
features
\begin{itemize}
\item $e:Y_{1,\beta}\rightarrow Y$ - the evaluation map. \item
$\psi$ - the first chern class of the cotangent line bundle on
$Y_{1,\beta}$. \item $\cf: Y_{1,\beta}\rightarrow Y_{0,\beta}$ - the
forgetful morphism.
\end{itemize}
\noindent Let $\hbar$ be a formal variable and \[ \label{eq: J}
J_{\beta}(Y):=e_*\left(\frac{[Y_{1,\beta}]}{\hbar(\hbar-\psi)}\right)=\sum_{k=0}^{\infty}
\frac{1}{\hbar^{2+k}}e_*(\psi^k\cap [Y_{1,\beta}]). \] The sum is
finite for dimension reasons. Let $p=\{p_1,p_2,...,p_k\}$ be a nef
basis of $H^2(Y,\QQ)$. For $t=(t_0,t_1,...,t_k)$ let
$$tp:=t_0+\sum_{i=1}^{k}t_ip_i.$$ The $\mathcal D$-module
for the quantum cohomology of $Y$ is generated by \cite{[G2]}
$$J(Y)=\exp\left(\frac{tp}{\hb}\right)\sum_{\beta\in
H_2(Y,\ZZ)}q^{\beta}J_{\beta}(Y)$$ where we use the convention
$J_0=1$. The generator $J(Y)$ encodes {\it all} of the one marking
Gromov-Witten invariants and gravitational descendants of $Y$.

\noindent For any ring $\mathcal A$, the formal completion of
$\mathcal A$ along the semigroup $MY$ of the rational curves of $Y$
is defined to be
\begin{eqnarray}
\mathcal A[[q^{\beta}]]:=\{\sum_{\beta\in
\mathrm{MY}}a_{\beta}q^{\beta}, & a_{\beta}\in \mathcal A, &
\beta-\text{effective}\}.
\end{eqnarray}
where $\beta\in H_2(Y,\ZZ)$ is {\it effective} if it is a positive
linear combination of rational curves. This new ring behaves like a
power series since for each $\beta$, the set of $\alpha$ such that
$\alpha$ and $\beta- \alpha$ are both effective is finite.
Alternatively, we may identify $q^{\beta}$ with $q_1^{d_1}\cdot
...\cdot q_k^{d_k}=\exp(t_1d_1+...+t_kd_k)$ where
$\{d_1,d_2,...,d_k\}$ are the coordinates of $\beta$ relative to the
dual basis of $\{p_1,...,p_k\}$.

\noindent We regard $J(Y)$ as an element of
$H^*(Y,\QQ)[t][[q^\beta]]$.

\noindent Let $*$ denote the small quantum product of $Y$. The small
quantum cohomology ring $(QH^*_sY,*)$ is a deformation of the
cohomology ring $(H^*(Y,\QQ[q^{\beta}]),\cup)$. Its structural
constants are three point Gromov-Witten invariants.

\noindent The generator $J(Y)$ arises naturally in the solution of
the small quantum differential equation: $$1\leq i\leq
k,~\displaystyle{\hbar
\partial /
\partial {t_i}=p_i
*}.$$ Furthermore $J(Y)$ may be used to generate
relations in $QH_s^*Y$. Let $\mathcal P(\hbar, \hbar \partial /
\partial t_i, q_i)$ be a polynomial differential operator where $q_i$ and $\hb$ act via
multiplication and $q_i=e^{t_i}$ are on the left of derivatives. If
$\mathcal P(\hbar, \hbar
\partial /
\partial t_i, q_i)J(Y)=0$ then $\mathcal P(0, p_i,
q_i)=0$ in $QH_s^*Y.$

\noindent If $Y$ is a toric variety, $J(Y)$ is related to an
explicit hypergeometric series $I(Y)$ via a change of variables
(\cite{[G1]}, \cite{[LLY3]}). Furthermore, if $Y$ is Fano then the
change of variables is trivial, i.e. $J(Y)=I(Y)$ thus completely
determining the one point Gromov-Witten invariants and gravitational
descendants of $Y$. In \cite{[E]} we conjectured an extension of
this result in the case of a projective bundle. We recall this
conjecture below.

\noindent Let $X$ be a projective manifold. Following Grothendieck's
notation, let $\pi: {\PP}(V)=\PP(\oplus_{j=0}^{n} L_j)\rightarrow X$
be the projective bundle of hyperplanes of a vector bundle $V$.
Assume that $L_0=\mathcal O_X$. The $H^*X$-module $H^*{\PP}(V)$ is
generated by $z:=c_1(\mathcal O_{\PP(V)}(1))$ with the relation
$$\prod_{i=0}^{n}(z-c_1(L_i))=0.$$ Let $s_i: X\rightarrow \PP(V)$ be the section
of $\pi$ determined by the $i$-th summand of $V$ and $X_i:=s_i(X)$.
Then $\mathcal O_{\PP(V)}(1)|_{X_i}\simeq L_i$. Let
$\{p_1,...,p_k\}$ be a nef basis of $H^2(X,\QQ)$. In \cite{[E]} we
showed that
\begin{lemma} If the line bundles $L_i, i=1,...,n$ are nef then

\noindent (a) $\{p_1,...,p_k,z\}$ is a nef basis of
$H^2(\PP(V),\QQ)$.

\noindent (b) The Mori cones of $X$ and $\PP(V)$ are related via
$$ M\PP(V)=MX\oplus \ZZ_{\geq 0}\cdot [\text{l}]$$
where $[\text{l}]$ is the class of a line in the fiber of $\pi$.
\end{lemma}

\noindent Here $MX$ is embedded in $M\PP(V)$ via the section $s_0$.
If $C\in \PP(V)$ is a rational curve, there exists a unique pair
$(\nu\geq 0,\beta\in MX$) such that $[C]=\nu[l]+\beta$. We will
identify the homology class $[C]$ with $(\nu,\beta)$. The generator
$J_{\PP(V)}$ is an element of
$H^*(\PP(V),\QQ)[t,t_{k+1}][[q_1^{\nu},q_2^\beta]].$

\noindent For a line bundle $L$ and a curve $\alpha$ we denote
$L(\alpha):=c_1(L)\cdot \alpha$.

\noindent Define the ``twisting'' factor
$$\cT_{\nu,\beta}:=\prod_{i=0}^n{\frac{\prod_{m=-\infty}^{0}(z-c_1(
L_i)+m\hb)}{ \prod_{m=-\infty}^{\nu-L_i(\beta)}(z-c_1(
L_i)+m{\hb})}}.$$ Let $I_{\nu,\beta}:=\cT_{\nu,\beta}\cdot
{\pi}^*J_{\beta}$ where $\pi^*$ is the flat pull back and define a
``twisted" hypergeometric series for the projective bundle $\PP(V)$:
$$I({\PP}(V)):=\exp\left(\frac{tp+t_{k+1}z}{\hb}\right)\cdot\sum_{\nu,\beta}{{q_1}^{\nu}}{{q_2}^{\beta}}I_{\nu,\beta}.$$

\noindent In \cite{[E]} we proposed the following
\begin{conjecture}
Let $L_i, ~i=1,...,n$ be nef line bundles such that $-K_X-c_1(V)$ is
ample. Then $J(\PP(V))=I(\PP(V)).$
\end{conjecture}

\noindent In the next section we prove this conjecture when $X$ is a
complete intersection in a toric variety.

\noindent In the last section we study the consequences of the
proposed conjecture in the relation between $QH^*_sX$ and
$QH^*_s\PP(V)$. Recall that $H^*{\PP}(V)$ is an $H^*X$-module
generated by $z:=c_1(\mathcal O_{\PP(V)}(1))$ with the relation
\begin{equation}
\label{eq: projbundlerel} \prod_{i=0}^{n}(z-c_1(L_i))=0.
\end{equation}
\noindent We show that the relations of $QH_s^*X$ that come from the
quantum differential equations lift to relations in $QH_s^*\PP(V)$.
We also show that (\ref{eq: projbundlerel}) deforms into the
relation
$$z \prod_{i=1}^{n}(z-c_1(L_i))=q_1$$ in $QH_s^*\PP(V)$.

\section{Toric case proof}

\noindent {\it \bf Toric varieties and torus actions.} Assume $Y$ is
a toric variety determined by a fan $\Sigma\subset \ZZ^m$. Denote by
$b_1,...,b_{r=m+k}$ its one dimensional cones. Let $Z_Y\subset
\CC^r$ be the variety whose ideal is generated by the products of
those variables which do {\it not} generate a cone in $\Sigma$. The
toric variety $Y$ is the geometric quotient of $\CC^r-Z(\Sigma)$ by
a torus of dimension $k$ \cite{[A]} \cite{[C]}.

\noindent Let $\tilde L_i,~i=0,1,...,n$ be toric line bundles and
$\tilde V=\oplus_{i=0}^n \tilde L_i$. The projective bundle $\pi:
\PP(\tilde V)\rightarrow Y$ is also a toric variety and there is a
canonical way to obtain its fan \cite{[O]}. Let $\ZZ^n$ be a new
lattice with basis $\{f_1,...,f_n\}$. The edges $b_1,...,b_r$ of
$\Sigma$ are lifted to new edges $B_1,B_2,...,B_r$ in $\ZZ^m\oplus
\ZZ^n$ and subsequently $\Sigma$ is lifted in a new fan $\Sigma_1$
in the obvious way. Let $\Sigma_2\subset 0\oplus \ZZ^n$ be the fan
of $\PP^n$ with edges $F_0=-\sum_{i=1}^n f_i,F_1=f_1,...,F_n=f_n$.
The canonical fan associated to $\PP(V)$ consists of the cones
$\sigma_1+ \sigma_2$ where $\sigma_1, \sigma_2$ are cones in
$\Sigma_1,\Sigma_2$. Let $N=r+n+1$. The torus $\TT=(\CC^*)^N$ acts
on both $Y$ and $\PP(V)$ by scaling of coordinates in respectively
$\CC^r$ and $\CC^N$. The one dimensional cones correspond to $\TT$
invariant divisors. The edges $F_i$ in the canonical fan of the
projective bundle, correspond to the divisors $z-c_1(L_i)$ where
$z:=c_1(\co_{\PP(V)}(1))$ while $B_i$ correspond to the pullback of
the base divisors associated with $b_i$ \cite{[M]}.

\noindent Let $L'_a: a=1,2,...,M$ be semiample line bundles (i.e.
generated by sections). Let $X$ be the zero locus of a generic
section s of $E=\oplus_{a=1}^{M}L'_a$ and let $V$ be the restriction
of $\tilde V$ to $X$. The total space of $\PP(V)$ is the zero locus
of the section $\pi^*(s)$ of the pull back bundle $\pi^*(\tilde V)$.
To assure that the conditions of the conjecture are met for the
bundle $\PP(V)$ over $X$ we assume that $\tilde L_i, i=1,1,...,n$
are nef and $-K_Y-\sum_{a=1}^M c_1(L'_a)-\sum_{i=0}^n c_1(\tilde
L_i)$ is ample.

\noindent Let $E_d$ be the bundle on $Y_{1,d}$ whose fiber over the
moduli point $(C,x_1,f:C\rightarrow Y)$ is $\oplus_a
H^0(f^*(L'_a))$. Denote by $s_E$ its canonical section induced by
$s$, i.e. $$s_E((C,x_1,f))=f^*(s).$$ The stack theoretic zero
section of $s_E$ is the disjoint union
\begin{equation}
\label{eq: zerolocus} Z(s_E)=\coprod_{i_*(\beta)=d}X_{1,\beta}.
\end{equation} The map $i_*: H_2X\rightarrow H_2Y$ is not injective in
general, hence the zero locus $Z(s_E)$ may have more then one
connected component. An example is the quadric surface in $\PP^3$.
The sum of the virtual fundamental classes $[X_{1,\beta}]$ is the
refined top Chern class of $E_d$ with respect to $s_E$.

\noindent There is a stack morphism $\PP(\tilde
V)_{1,(\nu,d)}\rightarrow Y_{1,d}.$ Let $\tilde E_{\nu,d}$ and
$\tilde s_E$ be the pull backs of $E_d$ and $s_E$. The zero section
of $\tilde s_E$ is the disjoint union $$z(\tilde
s_E)=\coprod_{i_*(\beta)=d} \PP(V)_{1,(\nu,\beta)}.$$ It follows
that
$$\sum_{i_*(\beta)=d}[\PP(V)_{1,(\nu,\beta)}]=c_{\text{top}}(\tilde E_{\nu,d})\cap
[\PP(\tilde V)_{1,(\nu,d)}].$$ Consider the following generating
functions
$$J(\PP(\tilde V),E)=\exp\left(\frac{tp+t_{k+1}z}{\hb}\right)
\sum{{q_1}^{\nu}}{{q_2}^d}e_*\left(\frac{c_{\text{top}}(\tilde
E_{\nu,d})\cap [\PP(\tilde V)_{1,(\nu,d)}]}{\hbar(\hbar-c)}\right)$$
and
$$\tilde I(\PP(\tilde V),E)=\exp\left(\frac{tp+t_{k+1}z}{\hb}\right)\sum q_1^{\nu}{{q_2}^d}~\cT_{\nu,d}~\pi^*e_*\left(\frac{c_{\text{top}}(
E_d)\cap [Y_{1,d}]}{\hbar(\hbar-c)}\right).$$

\begin{proposition} If $-K_Y-\sum_{a=1}^M
c_1(L'_a)-\sum_{i=0}^n c_1(\tilde L_i)$ is ample then
$$J(\PP(\tilde V),E)=\tilde I(\PP(\tilde V),E)$$
\end{proposition}
\begin{proof} Let $$I_d(Y,E)=\prod_a
\frac{\prod_{m=-\infty}^{L'_a(
d)}(L'_a+m\hb)}{\prod_{m=-\infty}^{0}(L'_a+m\hb)}\prod_i
\frac{\prod_{m=-\infty}^{0}(B_i+m\hb)}{\prod_{m=-\infty}^{B_i(
d)}(B_i+m\hb)}.$$From \cite{[G1]},\cite{[LLY2]},\cite{[LLY3]} we
know that $J(\PP(\tilde V),E)$ is related via a mirror
transformation to
$$I(\PP(\tilde V),E)=\exp\left(\frac{tp+t_{k+1}z}{\hb}\right)\cdot \sum {{q_1}^{\nu}}{{q_2}^d}\cT_{\nu,d}I_d(Y,E).$$ Likewise
$$J(Y,E)=\exp\left(\frac{tp}{\hb}\right)\sum
{{q_2}^d}e_*\left(\frac{c_{\text{top}}( E_d)\cap
[Y_{1,d}]}{\hbar(\hbar-c)}\right)$$ is related to
$$I(Y,E)=\exp\left(\frac{tp}{\hb}\right) \sum {{q_2}^d}I_d(Y,E).$$ Since $-K_{\PP(\tilde V)}-\sum_{a} c_1(L'_a)$
and $-K_Y-\sum_{a} c_1(L'_a)$ are ample, the mirror transformations
are particularly simple. Indeed, both series can be written as power
series of $\hb^{-1}$ as follows:
$$I(\PP(\tilde V),E)=1+\frac{P_1(q_1,q_2)}{\hb}+o(\hb^{-1}),~I(Y,E)=1+\frac{P_2(q_2)}{\hb}+o(\hb^{-1}),$$

\noindent where $P_1(q_1,q_2),P_2(q_2)$ are both polynomials
supported respectively in
$$\Lambda_1:=\{(\nu,d)~ | ~(-K_{\PP(\tilde V)}-\sum c_1(L'_a))=1;$$
$$z-c_1(\tilde L_j) \geq 0 ~\forall j=1,2,...,n; ~B_i
\geq 0~ \forall i=1,2,...,r\}$$ and
$$\Lambda_2:=\{d ~|~ (-K_{Y}-\sum c_1(L'_a))=1; B_i \geq 0~\forall i=1,2,...,r\}.$$
Then
$$J(\PP(\tilde V),E)=\exp\left(\frac{-P_1(q_1,q_2)}{\hb}\right)
I(\PP(\tilde V),E)$$ and
$$J(Y,E)=\exp\left(\frac{-P_2(q_2)}{\hb}\right)I(Y,E).$$

\noindent Let us examine the relation between $\Lambda_1$ and
$\Lambda_2$. From the exact sequence
\[0\rightarrow {\mathcal
O}_{\PP(V)}\rightarrow \pi^*V^*(1)\rightarrow T_{\PP(V)}\rightarrow
\pi^*T_X\rightarrow 0\] we find that
$$-K_{\PP(\tilde V)}-\sum c_1(L'_a)=-K_Y-\sum_i(\tilde
L_i)-\sum c_1(L'_a)+(n+1)z.$$ Assume $(\nu,d)\in \Lambda_1$. Since
$-K_Y-\sum_{a=1}^M c_1(L'_a)-\sum_{i=0}^n c_1(\tilde L_i)$ is ample
then $\nu=0$. Now $(z-c_1(\tilde L_j))\cdot (0,d)\geq 0~ \forall
j=1,2,...,n$ and $\tilde L_j$ are semipositive $\forall
j=1,2,...,n$. So $c_1(\tilde L_j)\cdot d=0~ \forall j=1,2,...,n$. It
follows that
$$(-K_{Y}-\sum c_1(L'_a))\cdot d=(-K_{\PP(\tilde V)}-\sum
c_1(L'_a))\cdot d=1,$$ so $d \in \Lambda_2$. Conversely, let $d \in
\Lambda_2$. Then $c_1(\tilde L_j)\cdot d=0~\forall j=1,2,...,n$
since $-K_Y-\sum_{a=1}^M c_1(L'_a)-\sum_{i=0}^n c_1(\tilde L_i)$ is
ample. It follows that
$$(-K_{\PP(\tilde V)}-\sum c_1(L'_a))\cdot d=(-K_{Y}-\sum c_1(L'_a))\cdot
d=1$$ and $$(z-c_1(\tilde L_j))\cdot d=0, \forall j=1,2,...,n$$ so
$(0,d)\in \Lambda_1$.

\noindent We have thus shown that $c_1(\tilde L_j)\cdot d=0, \forall
d \in \Lambda_2, \forall j=1,2,...,n$ and
$$\Lambda_1=\{(0,d)~|~d \in \Lambda_2\}.$$ \noindent It follows that
$\cT_{0,d}=1, \forall d\in \Lambda_2$ hence $P_1(q_1,q_2)=P_2(q_2)$.
Notice also that if we expand
$$\exp\left(\frac{-P_2(q_2)}{\hb}\right)=\sum_{\alpha}c_{\alpha}q_2^{\alpha}$$ then $$c_1(\tilde L_j)\cdot \alpha=0, \forall
j=1,2,...,n.$$ Hence for each $(\nu,d)\in M\PP(\tilde V)$ we have
$$\cT_{\nu,d}=\cT_{\nu,d+\alpha}.$$ Now the proposition follows easily.
\end{proof}

\noindent As we commented in the sentence after equation (\ref{eq:
zerolocus}), in general the above proposition is not very relevant
for our purpose. However, if we assume that the map
\begin{equation}
\label{eq: condition} i_*: H_2(X)\rightarrow H_2(Y)
\end{equation} is injective, then one can easily
show that
$$i_*(J_{\nu,\beta}(\PP(V)))=J_{\nu,i_*(\beta)}(\PP(\tilde V),E)$$ and
$$i_*(I_{\nu,\beta}(\PP(V)))= \tilde I_{\nu,i_*(\beta)}(\PP(\tilde
V),E).$$ The proposition shows that the conjecture holds for
complete intersection in toric varieties that satisfy condition
(\ref{eq: condition}).

\section{Relations in the small quantum cohomology ring}

\noindent In this section we use the proposed conjecture to study
small quantum deformations of the cohomological relation
$$H^*(\PP(V)=H^*X/\langle\prod_{i=0}^n (z-c_1(L_i))\rangle=0.$$ As explained in the
introduction, some of the relations in the small quantum cohomology
ring come from differential operators. Let
$$c_1(L_i)=\sum_{j=1}^k a_{ij}p_j, ~i=0,1,...,n.$$

\noindent We obtain two kinds of relations in $QH^*_s\PP(V)$.

\noindent First, the relations in $QH^*_sX$ that come from
differential operators may be lifted to relations in $QH^*_s\PP(V)$.
Indeed, consider a polynomial differential operator $$\mathcal
P(\hb, \hb
\partial /
\partial t_1,...,\hb
\partial /
\partial t_k, q_2)=\sum_{\alpha\in \Lambda}q_2^{\alpha}\mathcal P_{\alpha}$$ where $\Lambda\subset MX$ is a finite set.
Suppose that $$0=\mathcal P J(X)=\sum_{\alpha\in
\Lambda}q_2^{\alpha}\sum_{\beta}\mathcal
P_{\alpha}\left(\exp(\frac{pt}{\hb})q_2^{\beta}\right)J_{\beta}(X)$$
$$=\sum_{\alpha\in
\Lambda}q_2^{\alpha}\sum_{\beta}c_{\alpha,
\beta}\exp(\frac{pt}{\hb})q_2^{\beta}J_{\beta}(X)=\exp(\frac{pt}{\hb})\sum_{\alpha\in
\Lambda,\beta}q_2^{\alpha+\beta}c_{\alpha, \beta}J_{\beta}(X).$$ Let
$$\delta_{\alpha}=\prod_{i=1}^n\prod_{0}^{L_i\cdot \alpha-1}(\hb\frac{\partial}{\partial
t_{k+1}}-\sum_{j=1}^k a_{ij}\hb\frac{\partial}{\partial
t_j}-r_i\hb),~\pl=\sum_{\alpha\in
\Lambda}q_2^{\alpha}\delta_{\alpha}\mathcal P_{\alpha},$$ with the
convention that the factors corresponding to $L_i$ are missing if
$L_i(\alpha)=0$. We compute
$$\pl J(\PP(V))=\sum_{\alpha\in
\Lambda}q_2^{\alpha}\delta_{\alpha}\sum_{\nu,\beta}\mathcal
P_{\alpha}\left(q_2^{\beta}\exp(\frac{pt+zt_{k+1}}{\hb})\right)q_1^{\nu}\cT_{\nu,\beta}J_{\beta}=$$

$$\sum_{\alpha\in
\Lambda}q_2^{\alpha}\delta_{\alpha}\sum_{\nu,\beta}c_{\alpha,\beta}
\exp(\frac{pt+zt_{k+1}}{\hb})q_1^{\nu}q_2^{\beta}\cT_{\nu,\beta}J_{\beta}=0.$$
A simple calculation shows that
$$\delta_{\alpha}\left(\exp(\frac{pt+zt_{k+1}}{\hb})q_1^{\nu}q_2^{\beta}\cT_{\nu,\beta}\right)=\exp(\frac{pt+zt_{k+1}}{\hb})q_1^{\nu}q_2^{\beta}\cT_{\nu,\alpha+\beta}.$$
It follows that $$\pl
J(\PP(V))=\exp(\frac{pt+zt_{k+1}}{\hb})\sum_{\nu}q_1^{\nu}\sum_{\alpha\in
\Lambda,\beta}c_{\alpha,\beta}q_2^{\alpha+\beta}\cT_{\nu,\alpha+\beta}J_{\beta}(X)=0.$$
Hence the relation $\mathcal P(0, p_1,...,p_k, q_2)=0$ in $QH^*_sX$
lifts into the relation
$$\mathcal P(0, p_1,...,p_k, q_2\prod_{i=1}^n(z-c_1(L_i))=0$$ in
$QH^*_s\PP(V)$, where
$$\left(\prod_{i=1}^n(z-c_1(L_i)\right)^{\alpha}:=\prod_{i=1}^n(z-c_1(L_i))^{L_i(\alpha)}, \forall \alpha\in MX.$$

\noindent Second, we derive a $q_1$-deformation of the relation
$\prod_i (z-c_1(L_i))=0$. Consider the operator
$$\Delta(\hb\frac{\partial}{\partial t_1},...,\hb\frac{\partial}{\partial
t_k},\hb\frac{\partial}{\partial
t_{k+1}},q_1):=\prod_{i=0}^n(\hb\frac{\partial}{\partial
t_{k+1}}-\sum_{j=1}^k a_{ij}\hb\frac{\partial}{\partial t_j})-q_1.$$
It is easy to show that it satisfies $$\Delta J(\PP(V))=0.$$ It
follows that $\Delta(p_1,...,p_k,z,q_1)=0,$ in $QH^*_s\PP(V)$ i.e.
$$\prod_{i=0}^n (z-c_1(L_i))=q_1.$$ Much like $z^{n+1}=q$ is the deformation in $QH^*_s\PP^n$ of
$z^{n+1}=0$, the above relation is the deformation of $\prod_{i=0}^n
(z-c_1(L_i))=0.$

\vspace{+10 pt}

\noindent Department of Mathematics and Statistics

\noindent American University

\noindent 4400 Massachusetts Ave

\noindent Washington, DC 20016

\noindent aelezi@american.edu

\end{document}